\documentclass[12pt]{article}
\usepackage[T2A]{fontenc}
\usepackage[cp1251]{inputenc}
\usepackage{amsfonts}
\usepackage[active]{srcltx}
\usepackage{amsmath}
\usepackage{amssymb}
\usepackage{cite}
\newenvironment{sciabstract}{%
\begin{quote} \footnotesize} %\bf
{\end{quote}}
\newtheorem{theorem}{Theorem}%[section]

\newtheorem{lemma}{Lemma}

\newtheorem{remark}{Remark}
\usepackage{graphicx}
\title{\bf Economic numerical method of solving coefficient
inverse problem for 3D wave equation }
\author
{A.S. Leonov$^1$, A.B. Bakushinsky$^2$\\
\\
\normalsize{$^{1}$National Research Nuclear University 'MEPhI',}\\
\normalsize{Moscow, Russian Federation}\\
\normalsize{$^{2}$Federal Research Center 'Computer Science and Control'}\\
\normalsize{of Russian Academy of Sciences,}
\normalsize{Institute for Systems Analysis}\\
\normalsize{Moscow, Russian Federation}\\
}
\date{}

\begin{document}

\maketitle

\begin{sciabstract}
An inverse problem of acoustic sounding is under consideration
in a form of 3D inverse coefficient problem for wave equation.
Unknown coefficient is the local propagation velocity of
vibrations, which is associated with inhomogeneities of the
medium. We are looking for this coefficient, knowing special
time integrals of the scattered wave field. In the article, a
new linear 3D Fredholm integral equation of the first kind is
introduced, of which it is possible to find the unknown
coefficient from these time integrals. We present and
substantiate a numerical algorithm for solving this integral
equation. The algorithm does not require large computational
resources and big-time implementation. It is based on the use
of fast Fourier transform under some a priori assumptions about
unknown coefficient and observation region of the scattered
field. Typical results of solving this 3D inverse problem on a
personal computer for simulated data demonstrate the
capabilities of the proposed algorithm.

\textbf{Pacs}: 02.30.Zz, 02.30.Jr, 02.60.Cb

\end{sciabstract}

%\vspace{2mm}

\section{Introduction}

We study an inverse coefficient problems for the wave equation
\begin{equation} \label{EQ1}
\left\{\begin{array}{l} {\frac{1}{c^{2} (x)} u_{tt} (x,t)=
\Delta u(x,t)-g(t)\varphi (x),\, \, x\in \mathbb{R}^{3} ,\,
t>0} \\ \\{u(x,0)=u_{t} (x,0)=0,\, \, \, x=(x_{1} ,x_{2}
,x_{3} )\in \mathbb{R}^{3} } \end{array}\right.
\end{equation}
The initial value problem (\ref{EQ1}) can describe acoustic and
several electromagnetic wave processes. It is assumed below
that the coefficient $ c (x)> 0 $ is unknown in some bounded
region $ R\subset \mathbb{R}^{3} $. Without generality
restriction, we assume this area to be known. In the compliment
to the set $ R $, the value $c(x) $ is everywhere the same
constant $c_ {0}> 0 $, which is also known quantity.

The inverse coefficient problem for the equation (\ref{EQ1}) is
usually formulated as follows: knowing the function $ g (t), \,
\, \varphi(x) $, and the scattered field $u(x,t) $ in a domain
$Y\, (Y\subset \mathbb{R}^{3} ,\, R\cap Y=\emptyset )$, find
the function $c(x)$ in $R$. We assume further that the support
of the function $ \varphi(x) $ has no common points with the
sets $ Y $ and $ R $. In some cases, it is advisable to expand
this formulation of the inverse problem, believing that $ g
(t)$ and $\varphi (x) $ depend more on additional parameter $q
$. By doing so, we allow a plurality of experimental data
obtained from various types of sources.

Such inverse problem, in spite of the simplicity and
naturalness of its formulation, has been very difficult for
theoretical and numerical analysis. Various results relating to
this issue are reflected in a variety of articles and several
monographs. Having no possibility to consider in detail
previously obtained results, we note the most typical.

In the monograph \cite{1}, the original method of solving this
inverse problem proposed by M.V. Klibanov is described in
detail with its theoretical analysis. The impressive processing
results obtained with the help of this method for experimental
data are presented as well. However, apart from some specific a
priori assumptions on $c (x)$, the method requires significant
computing resources.

In the monographs \cite{2,3}, it is shown that universal
iterative algorithms studied there can be applied for solving
inverse problem in question. To effectively run such iterative
processes one must have a ''good'' initial guess. In addition,
a scant numerical implemen\-ta\-tion experience of these
algorithms in solving such inverse problem shows the need for
powerful computing system. Even if such a computer is
available, it is impossible to guarantee a reasonable time
calculations.

We note also the work \cite{4}, in which the inverse problem
for (\ref{EQ1}) is interpreted as a problem of optimal control,
and gradient methods for its solutions are formally applied.
Numerical results \cite{5} for the simulated data indicate the
following. If we use a sufficiently fine grid to approximate
the differential operators and have a supercomputer for
numerical realization of the methods, then having no
restrictions on the calculation time we can get a good
approximate solution. Thus, it is required powerful
computational resources in this case as well.

From the above brief review, it is evident that the inverse
problem for the equation (\ref{EQ1}) can not be considered
closed. It is still topical.

Quite a long time ago, it was observed that if one can observe
scattered field ''for a long enough time'' $t$ (formally, $ t
\in [0, \infty) $) and use specific time integrals of the
scattered field as the data to solve the inverse problem, one
can obtain a linear Fredholm integral equation of the first
kind (with the right-hand side including these integrals) to
find the unknown function $c(x)$. Such an equation is derived
in \cite{6}. Also, one can see there the references to previous
works. If we assume that we can register in the experiment the
mentioned integrals of the scattered field, not the field
itself, the inverse problem for (\ref{EQ1}) takes an entirely
different character. Firstly, the inverse problem becomes
linear. Secondly, the accumulation of information about the
scattered field in the form of these integrals allows us ''to
compress'' the data for solving the inverse problem, storage of
which requires substantial resources.

However, despite the linearity of the equation from \cite {6},
it numerical solution is also not an easy task because of the
need to solve ill-conditioned systems of linear  algebraic
equations (SLAE) having very high dimension.

The proposed work is adjacent to \cite{6}. In Sec.2, we obtain
a new linear integral equations of the first kind, from which
the function $c(x)$ can be found. Here we use some a priori
assumptions about the character of the inverse problem solution
$ c(x) $ together with certain assumptions about the sources of
the field $g(t)\varphi(x)$ and about properties of solutions to
the problem (\ref{EQ1}), $u(x,t)$. The right-hand sides of
these integral equations include the integral of the scattered
field $\int _{0}^{\infty }t^{2} u(x,t) dt$ or its partial
derivatives with respect to spatial variables. These values can
be registered experimentally. Besides, the right-hand side
contains some other integrals that can be computed ''a
priori''.

\begin{figure}[h]
  \centering
%%%%%%%%%%%%%%%%%%%%%%%%%%%%%%%%%%%%%%%%%%%%%
\includegraphics[width=100mm,height=100mm]{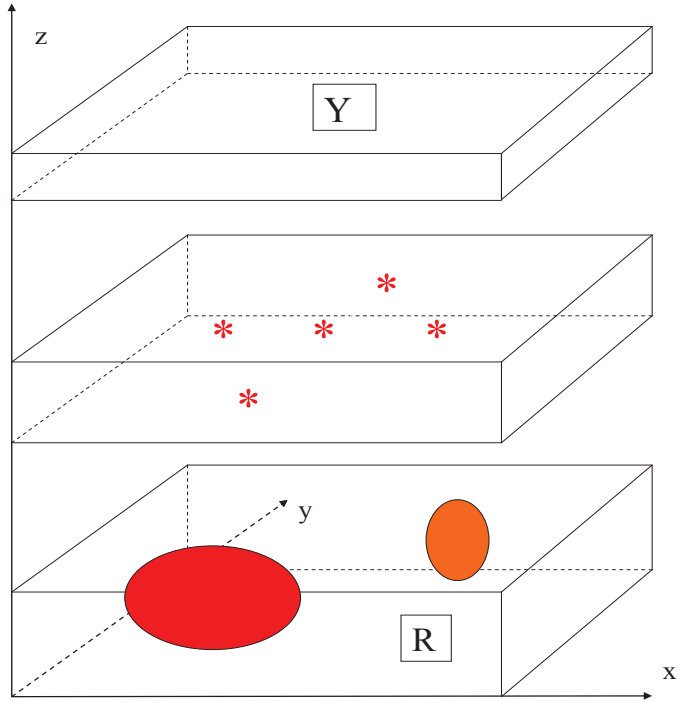}% bmp
%%%%%%%%%%%%%%%%%%%%%%%%%%%%%%%%%%%%%%%%%%%%%
  \caption{{\small Geometric registration scheme of the inverse problem data:
$R$ is a region of wave field scatterers, $Y$ is a domain of registration for the data
$V_{2} (x)$, «stars» are conditional positions of field sources.}}
  \label{fig1}
\end{figure}

Further, in Sec.3 we consider a special, but it seems to us,
accessible for implementing geometric registration scheme for
these integrals of the scattered field, so called
''registration in a flat layer'' (see Fig.\ref{fig1}).
Specifying the obtained integral equations for such a scheme,
we can see that one of them is more convenient to solve the
inverse problem, as the well-known methods for solving
ill-posed problems is easy to apply for it. It is this integral
equation we use further in the article. It can have more than
one solution as well as other resultant equation in Sec.3. In
this regard, we make the additional \emph{Assumption U}, which
actually defines a constraint on the unknown function $ c(x) $
and enables uniquely find $ c(x) $ from the integral equation
used.

In Sec.4, we present \emph{Algorithm 1} to solve such integral
equation. The algorithm is based on 2D Fourier transform of the
kernel and the right-hand side of the integral equation, so
that the original three-dimensional integral equation is
reduced to solving a set of one-dimensional integral equations
of the first kind by a selected method of regularization. Also,
we investigate the convergence of the proposed algorithm.

In Sec.5, we describe briefly a finite-dimensional
approximation of the considered integral equation and give a
finite-dimensional variant of \emph{Algorithm 1}. It is based
on the fast Fourier transform and on the reduction of the
three-dimensional discrete inverse problem to successive
solution of systems of linear algebraic equations. It turns out
that the application of this algorithm to determining the
coefficient $c(x)$ makes it possible to solve the problem
''quickly'' even by the use of a personal computer. Therefore,
we consider \emph{Algorithm 1} and its finite-dimensional
realization as the central result of the work. Finally, we
present in Sec. 6 a numerical illustration of our approach. In
particular, we estimate the efficiency of \emph{Algorithm 1}
analyzing its accuracy and operating speed for one of model
examples we have investigated.

\section{The integral equation for solving the inverse problem}

We suppose that the following \textbf{Assumptions} are
fulfilled.

1)The solution of (\ref{EQ1}) has the following smoothness:
$u(x,t)\in C^{2,2} (\mathbb{R}^{3} \times [0,\infty ))$.

2) The integrals
\[ V_{0} (x)=\int _{0}^{\infty }u(x,t)
dt,\, \,  V_{2} (x)=\int _{0}^{\infty }t^{2} u(x,t) dt,\]
converge for all $x\in \mathbb{R}^{3}$, and the functions
$V_{0}(x),~V_{2}(x)$ are regular at infinity ($\ |x|\rightarrow
\infty $) \cite[p.329]{TS}.

3) The equalities \[\int_{0}^{\infty }\Delta u(x,t)dt=\Delta
\left( \int_{0}^{\infty }u(x,t)dt\right) ,~\int_{0}^{\infty
}t^{2}\Delta u(x,t)dt=\Delta \left( \int_{0}^{\infty
}t^{2}u(x,t)dt\right); \] are valid.

4) $\mathop{\lim }\limits_{t\to \infty } \, \, t\,
u(x,t)=\mathop{\lim }\limits_{t\to \infty } t^{2} u_{t}
(x,\infty )=0,$ $\forall x\in \mathbb{R}^{3} $;

5) The function $\xi
(x)=\frac{1}{c_{0}^{2}}-\frac{1}{c^{2}(x)}$ is continuous in
$\mathbb{R}^{3}$ and compactly supported in $R$; the function
$\varphi (x)\in C^1(\mathbb{R}^{3} )$ is positive and compactly
supported in $D$; $g(t)\in C[0,+\infty )$ and the integrals
$A_0=\int _{0}^{\infty }g(t)dt$, $A_2=\int _{0}^{\infty }t^{2}
g(t)dt $ converge, and $A_0\neq 0$.

In fact, to meet the properties 1) - 4), the coefficients $
c(x)$, $\varphi (x)$ and $g(t) $ must satisfy more stringent
conditions than 5). However, a study of such conditions is a
separate scientific problem, which is actively investigated by
several authors (see, e.g. \cite {7,8,9}, and others). We are
not dealing with this issue, replacing its solution with the
requirements of 1) - 4). Note also that Assumptions 1) - 4) are
fulfilled for many functions $\varphi (x),\, \, g(t)$,
satisfying the conditions 5), if $c(x)=c_0$. For example, it is
true, if the function $ g (t) $ is finite or decreases
exponentially.

It follows from (\ref{EQ1}) and Assumptions 1) -- 3) and 5)
that
\[\frac{1}{c^{2} (x)} \int _{0}^{\infty }u_{tt} (x,t)dt =
\Delta \left(\int _{0}^{\infty }u(x,t) dt\right)-\varphi
(x)\int _{0} ^{\infty }g(t)dt .\] From here, integrating by parts the
member in the left side and taking into account the equality
\[\int _{0}^{\infty }u_{tt} (x,t)dt =u_{t} (x,\infty )-u_{t} (x,0)=0,\]
which follows from 1), 4), we obtain: $\Delta V_{0} (x)=A_{0}
\varphi (x),\, \, x\in \mathbb{R}^{3}$. It follows from
Assumption 2) that this Poisson equation has a classical
solution:
\[V_{0} (x)=-\frac{A_0}{4\pi } \int _{D }
\frac{\varphi (x')dx'}{\left|x-x'\right|}\in
C^2(\mathbb{R}^{3})  .\] Note that Assumption 5) entails that
$V_{0} (x)\neq 0,\,\forall x\in \mathbb{R}^{3}$.

Similarly, we can deduce from (\ref {EQ1}) and Assumptions 1) -
5) the relation
\[\frac{1}{c^{2} (x)} \int _{0}^{\infty }t^{2} u_{tt} (x,t)dt
=\Delta \left(\int _{0}^{\infty }t^{2} u(x,t) dt\right)-\varphi
(x) \int _{0}^{\infty }t^{2} g(t)dt. \] Together with the equality
\[\int _{0}^{\infty }t^{2} u_{tt} (x,t)dt =\left. t^{2}
u_{t} (x,t)\right|_{0}^{\infty } -2\int _{0}^{\infty }tu_{t}
(x,t)dt =-2\left. tu(x,t)\right|_{0}^{\infty } +2\int _{0}
^{\infty }u(x,t)dt =2V_{0} (x),\] it yields:
\[\frac{2V_{0} (x)}{c^{2} (x)} =\Delta V_{2} (x)-A_{2}
\varphi (x),\, \, x\in \mathbb{R}^{3} .\] Thus, introducing an
auxiliary function $\zeta (x)=\left(
\frac{1}{c_{0}^{2}}-\frac{1}{c^{2}(x)}\right) V_{0}(x)=\xi
(x)V_{0}(x)$, we obtain
\begin{equation}\label{eq33}
\zeta (x)-\frac{A_{2}}{2}\varphi (x)=\frac{1}{c_{0}^{2}}
V_{0}(x)-\frac{1}{2}\Delta V_{2}(x),\,\, x\in \mathbb{R}^{3}.
\end{equation}
The left-hand side of Equation (\ref{eq33}) is a function,
which under Assumption 5) is finite and integrable. So, there
exists its convolution with a locally integrable function
$\frac{1}{|x|}$ (see \cite[p.81]{VL}):
\begin{equation}\label{eq34}
\int_{\mathbb{R}^{3}}\left( \zeta (x^{\prime })-\frac{A_{2}}{2}\varphi
(x^{\prime })\right) \frac{dx^{\prime }}{|x-x^{\prime }|}=
\int_{\mathbb{R}^{3}}\left( \frac{1}{c_{0}^{2}}V_{0}(x^{\prime })-
\frac{1}{2}\Delta V_{2}(x^{\prime })\right) \frac{dx^{\prime
}}{|x-x^{\prime }|}.
\end{equation}
In addition, the regularity of the function $V_{2} (x)$ at
infinity (Assumption 2)), and third Green's formula \cite{TS}
lead to the equality
\begin{equation}\label{eq35}
\int_{\mathbb{R}^{3}}\frac{\Delta V_{2}(x^{\prime })d x^{\prime }}{|x-x^{\prime }|}=-4\pi V_{2}(x).
\end{equation}
From the relations (\ref{eq34}) and (\ref{eq35}), it follows
the existence of the convolution
$\int_{\mathbb{R}^{3}}\frac{V_{0}(x^{\prime })dx^{\prime
}}{|x-x^{\prime }|}$ and the validity of the equality
\begin{equation*} %\label{EQ2}
\int _{\mathbb{R}^3}\frac{\zeta (x')dx'}{\left|x-x'\right|}  =2\pi V_{2}
(x)+\frac{A_{2} }{2} \int _{\mathbb{R}^{3} }\frac{\varphi
(x')dx'}{\left|x-x'\right|} + \frac{1}{c_{0}^{2} }
\int _{\mathbb{R}^{3} }\frac{V_{0} (x')dx'}
{\left|x-x'\right|}  ,\, \, x\in \mathbb{R}^3.
\end{equation*}
If we know the function $V_{2} (x)$ in the domain $ Y $, we can
get by virtue of the finiteness of the functions $\zeta (x)$
and $\varphi (x)$ the following linear integral equation of the
first kind for the unknown $\zeta (x)=V_{0} (x)\xi (x)$, which
is associated with $c(x)$:
\begin{equation} \label{EQ2}
\int _{R}\frac{\zeta (x')dx'}{\left|x-x'\right|}  =2\pi V_{2}
(x)+\frac{A_{2} }{2} \int _{D }\frac{\varphi
(x')dx'}{\left|x-x'\right|} + \frac{1}{c_{0}^{2} }
\int _{\mathbb{R}^{3} }\frac{V_{0} (x')dx'}
{\left|x-x'\right|}  ,\, \, x\in Y,
\end{equation}
The integrals in the right-hand side of Equation (\ref{EQ2})
can be calculated because we know the values $\varphi
(x),A_{0},A_{2}$ a priori.

As a part of the inverse problem of acoustic sensing, the
measurements of the function $V_{2} (x)=\int _{0}^{\infty
}t^{2} u(x,t) dt,\, x\in Y,$ are associated with a special
accumulation of information about the sound pressure $u(x,t)$
at registration points, $x$. In principle, this can be done,
figuratively speaking, by processing the signals from a matrix
of microphones or other sensors located in $Y$ (analogously to
the registration of light signals by CCD camera). Note that
there is a family of so-called ''gradient microphones'', which
detect sound pressure gradient. With their help, it is possible
to measure the partial derivatives of the function $V_{2} (x)$
(for example, $\frac{\partial V_{2} (x)}{\partial x_{3} } $).
In this case, the other equation,
\begin{multline}\label{EQ3}
\int _{R}\frac{\partial }{\partial x_{3} }
\left(\frac{1}{\left|x-x'\right|} \right) \zeta (x')dx'=
 \\ =2\pi \frac{\partial V_{1} (x)}{\partial x_{3} }
  +\frac{A_2}{2} \int _{D }\frac{\partial }
  {\partial x_{3} } \left(\frac{1}{\left|x-x'\right|} \right)
  \varphi (x')dx'+ \frac{1}{c_{0}^{2} } \int _{\mathbb{R}^{3}
   }\frac{\partial }{\partial x_{3} } \left(\frac{1}
   {\left|x-x'\right|} \right)V_{0} (x')dx' ,\, \, x\in Y,
\end{multline}
can be obtained by analogy with (\ref{EQ2}). It differs from
Equation (\ref{EQ2}) in the kernel, that is used also at
calculating the integrals on the right. Equation (\ref{EQ3}) is
more convenient for solutions than (\ref{EQ2}). As shown below,
the well-known methods for solving ill-posed problems are
applicable for its solution. That is why we solve the equation
(\ref{EQ3}) in the sequel.
\begin{remark}
Instead of the functions $V_{2} (x)$, $\frac{\partial V_{2}
(x)}{\partial x_{3} } $, it is possible sometimes to measure
their analogs $V^{(0)}_{2} (x)$, $\frac{\partial V^{(0)}_{2}
(x)}{\partial x_{3} } $ for the case, when there are no
scatterers of the wave field in $R$. Formally, this corresponds
to the condition $\zeta(x)=0,\,\forall x\in R$. Then the
equations (\ref{EQ2}), (\ref{EQ3}) can be written as
\begin{multline} \label{EQ23}
\int _{R}\frac{\zeta (x')dx'}{\left|x-x'\right|}  =2\pi (V_{2}
(x)- V^{(0)}_{2} (x)),\\ \int _{R}\frac{\partial }{\partial
x_{3} } \left(\frac{1}{\left|x-x'\right|} \right) \zeta (x')dx'
=2\pi \left(\frac{\partial V_{2} (x)}{\partial x_{3} }-
\frac{\partial V^{(0)}_{2} (x)}{\partial x_{3} }\right),\,\
x\in Y.
\end{multline}
The integral equation (\ref{EQ23}) can be solved in the same
manner as discussed below in more detail Equation (\ref{EQ3}).
\end{remark}

\section{The scheme of data registration for the inverse problem
and special form of the basic integral equation}

The inverse problem of finding the function $c(x)$ is reduced
to solving a linear integral equation (\ref{EQ3}) with known
right-hand side, obtained from the measurement of the function
$\frac{\partial V_{2} (x)}{\partial x_{3} } $ in $ Y $ and
calculations of the integrals standing in right side of
(\ref{EQ3}). In this paper, we consider a specific scheme for
the registration  of the data, that is the function
$\frac{\partial V_{2} (x)}{\partial x_{3} } $, ''in a flat
layer''.

Below, for convenience we take for the variables $x_{1} ,x_{2}
,x_{3} $ the usual notation $x,y,z$. Figure 1 shows
schematically the geometry of the problem, with the region $R$
of heterogen\-ei\-ties, which scatter the incident waves, and
the domain $Y$, in which the scattered field is registered. The
domain $Y$ has a form of endless flat layer in variables $x,y$,
which is perpendicular to the axis $Oz$. The bounded region $R$
belongs to a similar layer. The asterisks shows conditionally
possible positions of the field sources.

Thus, we assume that $R=\mathbb{R}_{xy}^{2} \times [h,H],\, \,
Y=\mathbb{R}_{xy}^{2} \times [h_{Y} ,H_{Y} ]$ and consider
Equation (\ref{EQ3}) in the form
\begin{equation} \label{EQ4} \mathop{\int\!\!\!\!
\int\!\!\!\!\int}\nolimits _{R}K\left(x-x',y-y',z-z'\right)
\zeta (x',y',z')dx' dy'dz'=U(x,y,z),\, \, \, (x,y,z)\in Y,
\end{equation}
with the kernel \[K(x,y,z)=\frac{\partial }{\partial z}\left(
\frac{1}{\left( x^{2}+y^{2}+z^{2}\right) ^{1/2}}\right)
=-\frac{z}{\left( x^{2}+y^{2}+z^{2}\right) ^{3/2}}\] and the
right-hand side  $v$ given in (\ref{EQ3}). We transform
Equation (\ref{EQ4}) in the following way:
\[\int _{h}^{H}dz'\mathop{\int\!\!\!\!\int}\nolimits
_{\mathbb{R}_{xy}^{2} }K(x-x',y-y',z-z')\zeta (x',y',z')dx'dy'
=U(x,y,z),\, \, \,   (x,y)\in\mathbb{R}_{xy}^{2} ,\, z\in
[h_{Y} ,H_{Y} ].\]

Assumptions about the function $\zeta(x,y,z)$ and definition of
the functions $ K(x, y, z)$ and $U (x, y, z)$ ensure the
inclusions: $\zeta\in L_2(R),\,\,v\in L_2(Y)$;
$\zeta(x',y',z')\in L_2(\mathbb{R}_{xy}^{2}), \,\,v(x,y,z)\in
L_2(\mathbb{R}_{xy}^{2})$ for all admissible $z,z'$, as well as
the following relations: $||K(x-x',y-y',z-z')||_{L_2(Y\times
R)}<\infty$,
$||K(x,y,z-z')||_{L_2(\mathbb{R}_{xy}^{2})}<\infty$, $\forall
z,z',\,z\in[h_Y,H_Y],\,z'\in[h,H]$. Then, using the
two-dimensional Fourier transforms $\tilde{K}(\omega _{1}
,\omega _{2} ,z),\, \, \tilde{\zeta }(\omega _{1} ,\omega _{2}
,z)$, $\tilde{v}(\omega _{1} ,\omega _{2} ,z)$ of these
functions in the variables $(x,y)$ and taking into account the
convolution theorem, we obtain a family of one-dimensional
integral equations of the first kind
\begin{equation} \label{EQ5}
\int _{h}^{H}\tilde{K}(\omega _{1} ,\omega _{2} ,z-z')
\tilde{\zeta }(\omega _{1} ,\omega _{2} ,z')dz'=\tilde{v}
(\omega _{1} ,\omega _{2} ,z) ,\, \, z\in [h_{Y} ,H_{Y} ].
\end{equation}

The equation (\ref{EQ4}) can be written in the operator form $
A\zeta = v $ with linear and bounded integral operator $ A: L_2
(R)\rightarrow L_2(Y) $. Also, we represent Equations
(\ref{EQ5}) in the operator form $A(\omega _{1},\omega
_{2})\tilde{\zeta }=\tilde{v}(\omega _{1},\omega _{2},z)$ with
bounded linear integral operators $A(\omega _{1},\omega
_{2})[\cdot]=\int _{h}^{H}\tilde{K}(\omega _{1} ,\omega _{2}
,z-z')[\cdot]dz'$ acting from $L_2[h,H]$ to $L_2[h_Y,H_Y]$.

It is well known that the solution to the integral equation of
the form
\begin{equation}\label{EQ41}
\int _{R}\frac{\partial }{\partial x_{3} }
\left(\frac{1}{\left|x-x'\right|} \right) \zeta (x')dx' = v(x),\, \, x\in Y,
\end{equation}
may be non-unique in the class of functions $\zeta(x)$, which
we consider. In particular, this is true for the equation
(\ref{EQ4}). This means that the null-space $N(A)$ of the
operator $A$ contains non-trivial elements $\zeta=\zeta_0(x)$:
$A\zeta_0=0$.

\noindent \textbf{Assumption U.} We assume that the required
solution $\zeta=\bar\zeta(x)$ of the equation (\ref{EQ4})
satisfies the condition
$(\bar\zeta,\zeta_0)_{L_2(R)}=0,\,\forall \zeta_0\in N(A)$.

This is consistent with our desire to be free as possible in
the desired solution of the artefacts $\zeta_0(x)$, creating
the observed zero function $v(x)$. Thus, we seek the unique
normal solution $\zeta=\bar\zeta(x)$ of the linear operator
equation $A\zeta = v$, that is Equation (\ref{EQ4}), in the
space $L_2(R)$ assuming that this solution is continuous and
finite, $\mathrm{supp}\, \bar\zeta\subset R$. Also, we do not
exclude the case when $\bar\zeta(x)$ belongs to one of the
uniqueness classes of solutions to the equation (\ref{EQ41}),
which have been studied previously, for example, in the works
\cite{Str,Pril} and etc.

Currently, many stable methods, so called RAs, are known for
finding normal solutions $\zeta=\bar\zeta\in Z$ of linear
operator equations $A\zeta=v, \,v\in U,$ in Hilbert spaces
$Z,U$ (see, e.g. \cite{BG,Vai,Mor,EHN,TGSY} and etc). Suppose
that a parametric family of operators $R_\alpha (A):
U\rightarrow Z$ presents one of such methods. Assume that
instead of $v$, we have at our disposal its approximation, that
is an element $v_\delta\in U$, which meets the approximation
condition $||v-v_\delta||_{U}\le \delta$. Then, with a suitable
choice of the parameter $\alpha=\alpha(\delta)$ the convergence
of approximate solutions $\zeta_\delta=R_{\alpha(\delta)}
(A)v_\delta$ takes place:
$||\zeta_\delta-\bar\zeta||_{Z}\rightarrow 0$ as
$\delta\rightarrow 0$. In particular, the regularity conditions
for elements $\zeta_\delta$, that is
\begin{equation}\label{reg}
\overline{\lim_{\delta\rightarrow 0}}||\zeta_\delta||_Z\le||\bar\zeta||_{Z},\,\,
\lim_{\delta\rightarrow 0}||A\zeta_\delta-v_\delta||_U=0,
\end{equation}
are sufficient for such convergence (see, e.g. \cite{Mor,TLY}).

\section{Algorithm for solving the basic equation (\ref{EQ4}) and the stability of the algorithm }

\noindent\textbf{\emph{Algorithm 1}}.

\noindent 1)Calculation of 2D Fourier transforms
$\tilde{K}(\omega _{1} ,\omega _{2} ,z)$ и $\tilde{v}(\omega
_{1} ,\omega _{2} ,z)$ for each $z\in[h,H]$.

\noindent 2) For each $(\omega _{1} ,\omega _{2})$, finding an
approximate normal solution $\tilde{\zeta}_\delta(\omega _{1}
,\omega _{2} ,z)$ of the integral equation (\ref{EQ5}) via an
RA, which ensures the regularity conditions of the form
\begin{gather}
||\tilde{\zeta}_\delta(\omega_{1} ,\omega _{2}
,z)||_{L_2(R_{\omega_{1} \omega _{2}}^2)} \le
||\bar{\zeta}(\omega_{1} ,\omega _{2}
,z)||_{L_2(\mathbb{R}_{\omega_{1} \omega _{2}}^2)},\label{reg1} \\
\lim_{\delta \rightarrow 0}\left\Vert A(\omega _{1},\omega
_{2})\widetilde{\zeta }_{\delta }(\omega _{1},\omega
_{2},z)-\widetilde{v}_{\delta }(\omega _{1},\omega
_{2},z)\right\Vert _{L_{2}(\mathbb{R}_{\omega _{1}\omega
_{2}}^{2})}=0,~\ \forall z\in \lbrack h,H].\label{reg2}
\end{gather}

\noindent 3) For each $z\in[h,H]$, calculating an approximate
solution of the equation (\ref{EQ4}) using 2D inverse Fourier
transform $F^{-1}$ in the variables $(x,y)$:
\[\zeta _{\delta }(x,y,z)=F^{-1}[\tilde{\zeta }_{\delta }(\omega
_{1},\omega _{2},z)](x,y).\]

\noindent 4) Finding the function $\xi_\delta
(x,y,z)=\zeta_{\delta }(x,y,z)/V_{0} (x,y,z)$, which
approximates the value $\xi
(x,y,z)=\frac{1}{c_{0}^{2}}-\frac{1}{c^{2}(x,y,z)}$ and
calculating the approximation for $c(x,y,z)$ from the last
equality.

Note that the conditions (\ref{reg1}) and (\ref{reg2}) hold
true for Tikhonov regularization, if we select the parameter
$\alpha(\delta)$ using the discrepancy principle \cite{Mor} or
its generalization \cite{TLY}. The same property has the TSVD
method (see, e.g. \cite{Hans,EHN}) under suitable choice of the
regularization parameter.

Now, we turn to the justification of the stability of the
proposed \emph{Algorithm 1}.
\begin{lemma}\label{L1}
Let $\{\tilde{\zeta}_{\omega _{1}\omega _{2}}^{(norm)}(z)\}\in
L_{2}[h,H]$ be a family of normal solutions to the equations
(\ref{EQ5}) for all considered values $\omega _{1},\omega
_{2}$. Then the equality holds:
$\bar{\zeta}(x,y,z)=F^{-1}\left[ \tilde{\zeta}_{\omega
_{1}\omega _{2}}^{(norm)}(z)\right] (x,y)$.
\end{lemma}
\emph{Proof}. We introduce the element $\zeta
^{(norm)}(x,y,z)=F^{-1}\left[ \tilde{\zeta}_{\omega _{1}\omega
_{2}}^{(norm)}(z)\right] (x,y)$ and denote 2D Fourier transform
of the normal solution $\bar\zeta(x,y,z)$ to the equation
(\ref{EQ4}) as  $\tilde{\bar\zeta}(\omega _{1},\omega _{2},z)$.
The following Plancherel equalities are valid for all
$z\in[h,H]$:
\begin{multline*}
\iint\limits_{\mathbb{R}_{xy}^{2}}\left[ \zeta
^{(norm)}(x,y,z)\right]^{2} dxdy=\frac{1}{4\pi
^{2}}\iint\limits_{\mathbb{R}_{\omega _{1} \omega
_{2}}^{2}}\left[ \tilde{\zeta}_{\omega _{1}
\omega _{2}}^{(norm)}(z)\right] ^{2}d\omega _{1}d\omega _{2},\\
\iint\limits_{\mathbb{R}_{xy}^{2}}\left[
\bar{\zeta}(x,y,z)\right] ^{2} dxdy=\frac{1}{4\pi
^{2}}\iint\limits_{\mathbb{R}_{\omega _{1} \omega
_{2}}^{2}}\left[ \tilde{\bar\zeta}(\omega _{1},\omega _{2},z)
\right] ^{2}d\omega _{1}d\omega _{2}.
\end{multline*}
The function $\tilde{\bar\zeta}(\omega _{1},\omega _{2},z)$
satisfies the equations (\ref{EQ5}) for all $z\in[h,H]$.
Therefore
\[\left\Vert \zeta _{\omega _{1}\omega
_{2}}^{(norm)}(z)\right\Vert
_{L_{2}[h,H]}^{2}=\int\limits_{h}^{H}\left[
\tilde{\zeta}_{\omega _{1}\omega _{2}}^{(norm)}(z)\right]
^{2}dz\leq \int\limits_{h}^{H}\tilde{\bar\zeta}^{2}(\omega
_{1},\omega _{2},z)dz=\left\Vert \tilde{\bar\zeta}(\omega
_{1},\omega _{2},z)\right\Vert _{L_{2}[h,H]}^{2}.\] Combining
these relations and applying Fubini's theorem, we obtain
\begin{multline*}
\left\Vert \zeta ^{(norm)}(x,y,z)\right\Vert _{L_{2}(R)}^{2}
=\int\limits_{h}^{H}dz\iint\limits_{\mathbb{R}_{xy}^{2}}\left[
\zeta
^{(norm)}(x,y,z)\right] ^{2}dxdy=\\=\frac{1}{4\pi ^{2}}\int\limits_{h}^{H}dz%
\iint\limits_{\mathbb{R}_{\omega _{1}\omega _{2}}^{2}}\left[
\tilde{\zeta}_{\omega _{1}\omega _{2}}^{(norm)}(z)\right]
^{2}d\omega _{1}d\omega _{2}=\frac{1}{4\pi
^{2}}\iint\limits_{\mathbb{R}_{\omega _{1}\omega
_{2}}^{2}}\left\Vert \zeta _{\omega _{1}\omega
_{2}}^{(norm)}(z)\right\Vert
_{L_{2}[h,H]}^{2}d\omega _{1}d\omega _{2}\leq \\ \le \frac{1}{4\pi ^{2}}%
\iint\limits_{\mathbb{R}_{\omega _{1}\omega
_{2}}^{2}}\left\Vert \tilde{\bar\zeta}(\omega _{1},\omega
_{2},z)\right\Vert _{L_{2}[h,H]}^{2}d\omega _{1}d\omega _{2}=
\frac{1}{4\pi
^{2}}\int\limits_{h}^{H}dz\iint\limits_{\mathbb{R}_{\omega
_{1}\omega _{2}}^{2}}\left[ \tilde{\bar\zeta}(\omega
_{1},\omega _{2},z)\right]
^{2}d\omega _{1}d\omega _{2}= \\ =\int\limits_{h}^{H}dz\iint\limits_{\mathbb{R}_{xy}^{2}}%
\left[ \bar{\zeta}(x,y,z)\right] ^{2}dxdy=\left\Vert \bar{\zeta}%
(x,y,z)\right\Vert _{L_{2}(R)}^{2}.
\end{multline*}
By the uniqueness of the normal solution to the equation
(\ref{EQ4}), this implies equality $\bar{\zeta}(x,y,z)=\zeta
^{(norm)}(x,y,z)=F^{-1}\left[ \tilde{\zeta}_{\omega _{1}\omega
_{2}}^{(norm)}(z)\right] (x,y)$. $\square$
\begin{theorem}\label{TH1}
Algorithm 1 ensures the convergence $\left\Vert \zeta _{\delta
}(x,y,z)-\bar{\zeta}(x,y,z)\right\Vert _{L_{2}(R)}\rightarrow
0$ as $\delta \rightarrow 0$.
\end{theorem}
\emph{Proof}. We define for $z\in \lbrack h,H]$ the family of
the functions \[\eta _{\delta }(z)=\left\Vert
\tilde{\zeta}_{\delta }(\omega _{1},\omega
_{2},z)-\tilde{\zeta}_{\omega _{1}\omega
_{2}}^{(norm)}(z)\right\Vert _{L_{2}\left( \mathbb{R}_{\omega
_{1}\omega _{2}}^{2}\right) }^{2}=\left\Vert
\tilde{\zeta}_{\delta }(\omega _{1},\omega
_{2},z)-\tilde{\bar\zeta}(\omega _{1},\omega _{2},z)\right\Vert
_{L_{2}\left( \mathbb{R}_{\omega _{1}\omega _{2}}^{2}\right)
}^{2}.\] The dual form of this family's representation follows
from Lemma 1. The properties (\ref{reg1}) and (\ref{reg2}) of
the used RA guarantee for every $z\in[h,H]$ the convergence of
approximate solutions $\tilde{\zeta}_{\delta }(\omega
_{1},\omega _{2},z)$ to normal solutions $\tilde{\zeta}_{\omega
_{1}\omega _{2}}^{(norm)}(z)$ in the space $L_{2}\left(
\mathbb{R}_{\omega _{1}\omega _{2}}^{2}\right) $ as $\delta
\rightarrow 0$. Therefore, $\lim_{\delta \rightarrow 0}\eta
_{\delta }(z)=0,~\forall z\in \lbrack h,H]$. Besides, the
conditions (\ref{reg1}) and the inequality
\begin{multline*}
0\leq \eta _{\delta }(z)\leq \left( \left\Vert
\tilde{\zeta}_{\delta }(\omega _{1},\omega _{2},z)\right\Vert
_{L_{2}\left( \mathbb{R}_{\omega _{1}\omega _{2}}^{2}\right)
}+\left\Vert \tilde{\bar\zeta}(\omega _{1},\omega
_{2},z)\right\Vert _{L_{2}\left( \mathbb{R}_{\omega _{1}\omega
_{2}}^{2}\right) }\right) ^{2}\leq \\ \le 4\left\Vert
\tilde{\bar\zeta}(\omega _{1},\omega _{2},z)\right\Vert
_{L_{2}\left( \mathbb{R}_{\omega _{1}\omega _{2}}^{2}\right)
}^{2}\equiv s(z)
\end{multline*}
imply that the functions $\eta _{\delta }(z)$ have integrable
majorant $s(z)$. Then the convergence to be proved can be
derived from the corresponding Plancherel equality and
Lebesgue's theorem on passage to the limit in the following
way:
\begin{eqnarray*}
\lim_{\delta \rightarrow 0}\left\Vert \zeta _{\delta }(x,y,z)-\bar{\zeta}%
(x,y,z)\right\Vert _{L_{2}(R)}^{2} &=&\lim_{\delta \rightarrow 0}\frac{1}{%
4\pi ^{2}}\int\limits_{h}^{H}\left\Vert \tilde{\zeta}_{\delta }(\omega
_{1},\omega _{2},z)-\tilde{\bar\zeta}(\omega _{1},\omega _{2},z)
\right\Vert _{L_{2}\left( \mathbb{R}_{\omega _{1}\omega
_{2}}^{2}\right) }^{2}dz= \\
&=&\lim_{\delta \rightarrow 0}\frac{1}{4\pi ^{2}}\int\limits_{h}^{H}\eta
_{\delta }(z)dz=0. \,\,\square
\end{eqnarray*}

\section{Finite-dimensional approximation of the problem and
the numerical implementation of the algorithm}

We replace in (\ref{EQ3}) and (\ref{EQ4}) the space
$\mathbb{R}^{3} $ by the region $\Pi =[-r,r]\times [-r,r]\times
[-r,r]$ and the space $\mathbb{R}^{2}_{xy} $ by the rectangle
$\Pi_{xy} =[-r,r]\times [-r,r]$ with $r>0$ ''large enough''. We
carry out an approximation of equations (\ref{EQ4}) and
(\ref{EQ5}) in the domain $\Pi $ by the finite difference
method introducing uniform grids for $x,\omega _{1} \in [-r,r]$
and $y,\omega _{2} \in [-r,r]$ of the size $N$, as well as the
grids for $z,z'$: $\left\{z_{i} \right\}\in [h_{Y} ,H_{Y} ],\,
\, \left\{z'_{j} \right\}\in [h,H]$ of the size $M$ and $M'$
respectively. After that, we apply 2D fast Fourier transform in
the first section of \emph{Algorithm 1} for calculation of
discrete analogues of the functions $\tilde{K}(\omega _{1}
,\omega _{2} ,z)$, $\tilde{v}(\omega _{1} ,\omega _{2} ,z)$. In
the second section of \emph{Algorithm 1} we approximate
integrals in (\ref{EQ5}) using quadrature and obtain $N^{2}$
systems of linear algebraic equations for subsequent solution:
\begin{equation} \label{EQ6}
A^{(m)} \tilde{\zeta }^{(m)} =\tilde{U}^{(m)} ,\, \,
m=1,...,N^{2}.
\end{equation}
Here $A^{(m)} =\left[\nu _{ij} \tilde{K}(\omega _{1}^{(m)}
,\omega _{2}^{(m)} ,z_{i} -z'_{j} )\right]$ are matrices of the
size $M\times M'$ and $\tilde{U}^{(m)} =[\tilde{U}(\omega
_{1}^{(m)} ,\omega _{2}^{(m)} ,z_{i} )] $ are columns of the
length $M$. The values $\left(\omega _{1}^{(m)} ,\omega
_{2}^{(m)} \right)$ are grid points for the variables $(\omega
_{1} ,\omega _{2} )$ numbered by a single superscript $m$ and
$\nu _{ij}$ are quadrature coefficients. The SLAEs (\ref{EQ6})
was solved by application of the RAs with properties
(\ref{reg1}) and (\ref{reg2}).

In doing so, we used a number of RAs, namely, Tikhonov
regularization in the standard and iterated form, the TSVD
method and some others (for their implementation, see e.g.
\cite{2,3,TGSY,10,Hans,EHN,TLY}). Numerical experiments have
shown that the best results in the regularization of systems
(\ref{EQ6}) gives the TSVD method. The results of its
application are presented in the following section.

\section{Model example}

We write the equation (\ref{EQ3}) formally, regardless of
Assumptions 1) - 5), for the function
%---------------------------------------------------------
$\varphi (x)\, =\sum_{l=1}^L\delta (x-x_{l} ),\, \, x_{l}
\notin R,\, \, x_{l} \notin Y$.
%---------------------------------------------------------
Then
\[V_{0} (x)=-\frac{A_{0} }{4\pi }\sum_{l=1}^L {\frac{1}{\left|x-x_{l} \right|}} ,\]
and the equation (\ref{EQ3}) takes the form
\begin{multline}\label{EQ7}
\int_{R}\frac{\partial }{\partial x_{3}}\left(
\frac{1}{\left\vert
x-x^{\prime }\right\vert }\right) \zeta (x^{\prime })dx^{\prime }=\\=2\pi \frac{%
\partial V_{2}(x)}{\partial x_{3}}+\frac{A_{2}}{2}\sum_{l=1}^{L}\frac{%
\partial }{\partial x_{3}}\left( \frac{1}{\left\vert x-x_{l}\right\vert }%
\right) -\frac{A_{0}}{4\pi c_{0}^{2}}\sum_{l=1}^{L}\int_{\mathbb{R}^{3}}%
\frac{\partial }{\partial x_{3}}\left( \frac{1}{\left\vert
x-x^{\prime }\right\vert }\right) \frac{dx^{\prime
}}{\left\vert x^{\prime }-x_{l}\right\vert },\,\,\,x\in Y.
\end{multline}
This form of the integral equation can be used for the
finite-dimensional approximation as in Section 5, with
''sufficiently fine'' grids, meaning $ \delta (x-x_{l}) $ to be
some smooth $\delta $-shaped family of finite functions on
these grids

Keeping in mind the experimental scheme shown in Figure 1, we
define a model solution
\[\begin{array}{l} {\xi (x,y,z)=a_{1} \exp \left(-x^{2} -2y^{2}
\right){\rm +}a_{2} \exp \left[-3(x+4)^{2} -(y-5)^{2} +(x+4)(y-5)
\right]+} \\ {+a_{3} \exp \left\{-0.9\left[(x-4)^{2} -(y+4)^{2}
 +(x-4)(y+4)\right]\right\}} \end{array}\]
with
\[\begin{array}{l} { a_{1} =\left(1-4(z-1.5)^{2} \right)^{2} ,\, \,
a_{2} =0.4\max \left\{1-(z_{0} -1.3)^{2} {\rm ,\; }0\right\}},
\\ {a_{3} =0.2\max \left\{\left[1-(z_{0} -1.7)^{2} \right]^{{\rm 2}}
 {\rm ,\; }0\right\} } \end{array}\]
and $(x,y,z)\in R=[-10,10]\times [-10,10]\times [1,2]$.
Registration domain of the scattered wave field  is represented
as $Y=[-10,10]\times [-10,10]\times [6,7]$. We use ten $\delta
$-shaped sources whose positions are given by the points
\[\, (x_{l} ,y_{l} ,z_{l} )=(0,0,z_{0} ),\, (-r,0,z_{0} ),
(r,0,z_{0} ),(0,-r,z_{0} ),(0,r,z_{0} );\, \, \, z_{0} =3,5; \,
\, \, r=8.\]
\begin{figure}%[h]
  \centering
%%%%%%%%%%%%%%%%%%%%%%%%%%%%%%%%%%%%%%%%%%%%%
\includegraphics[width=35mm,height=55mm]{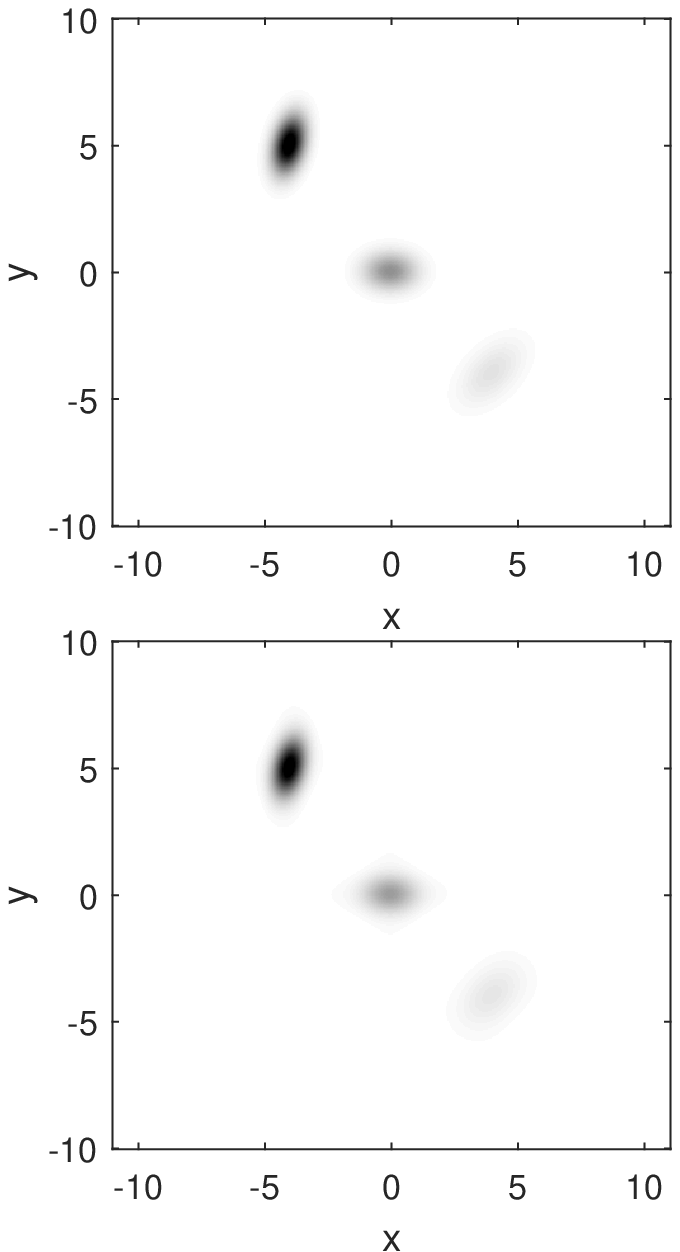}% bmp
\includegraphics[width=35mm,height=55mm]{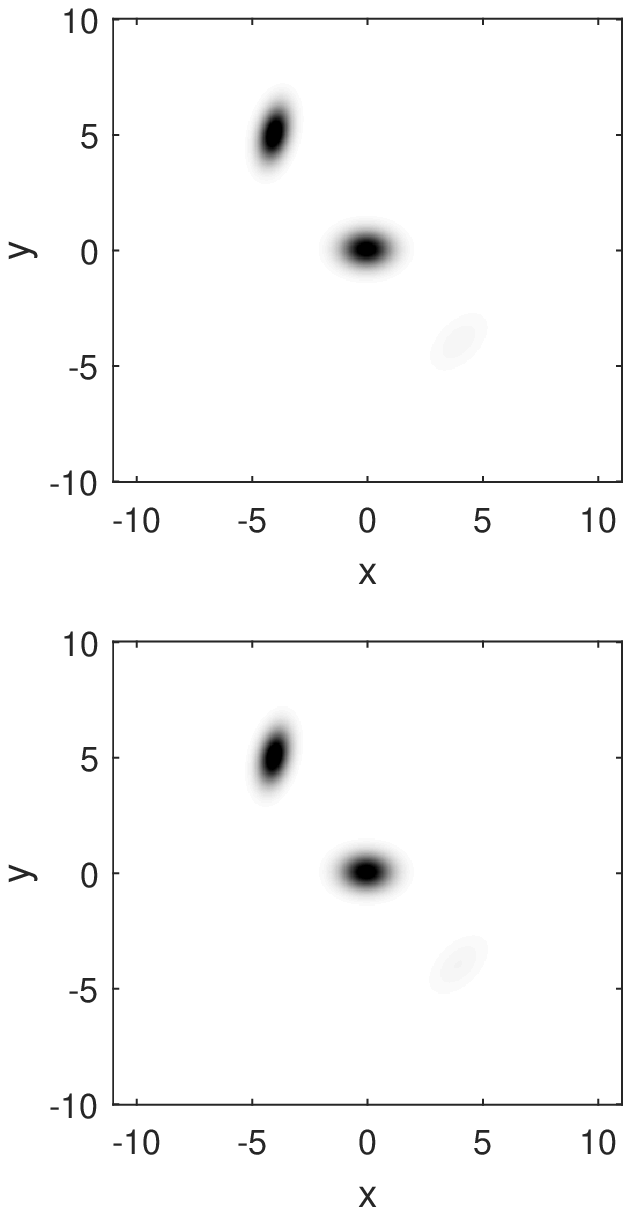}
\includegraphics[width=35mm,height=55mm]{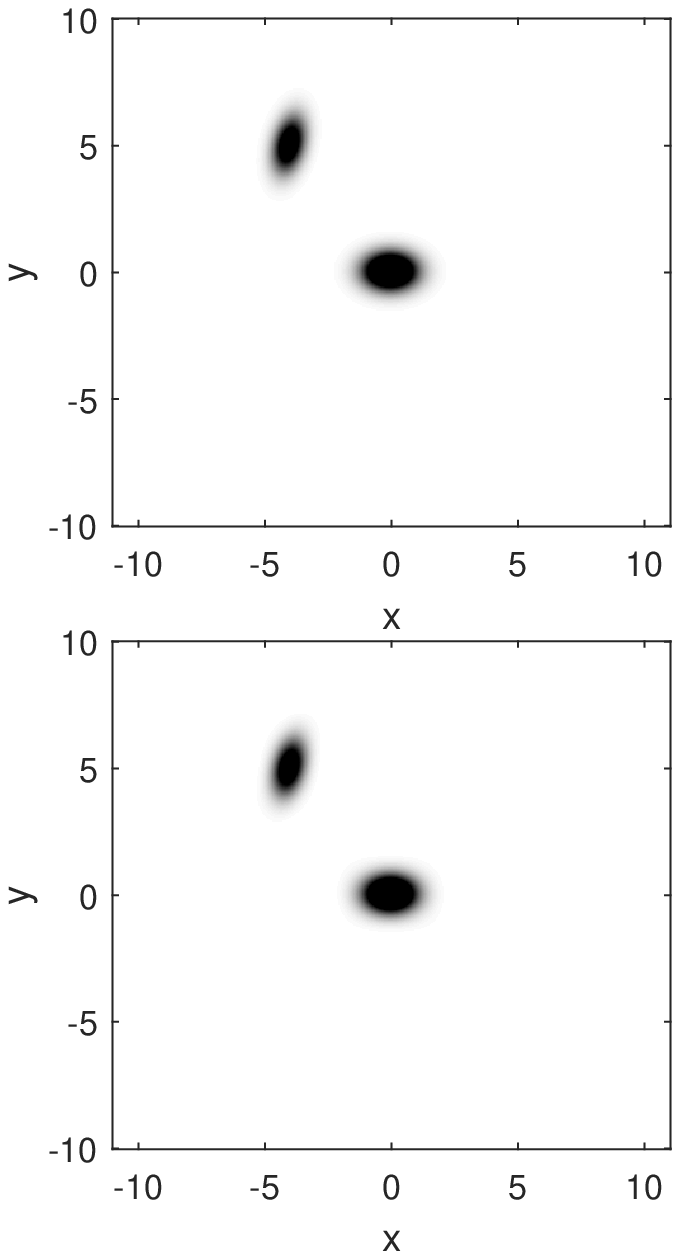}
\includegraphics[width=35mm,height=55mm]{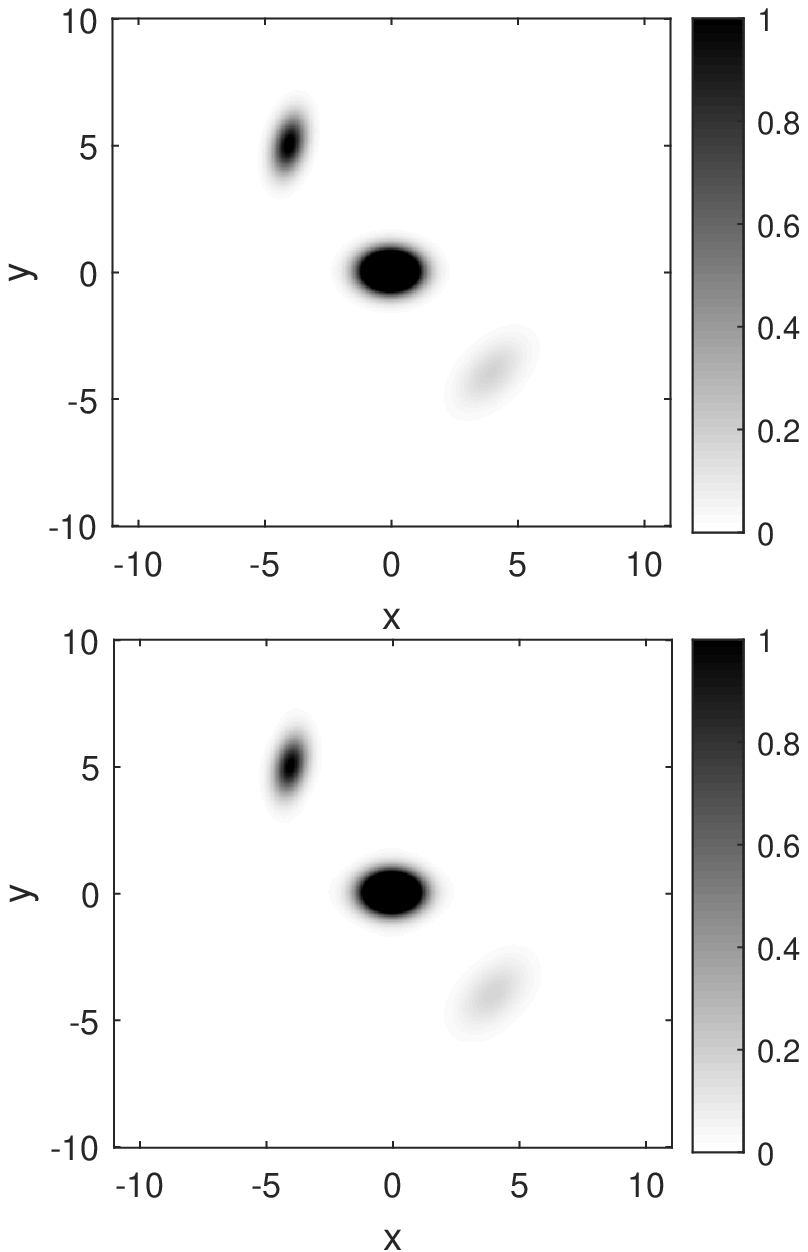}
%%%%%%%%%%%%%%%%%%%%%%%%%%%%%%%%%%%%%%%%%%%%%
  \caption{{\small A qualitative comparison of the exact (top) and an approximate
(bottom) solutions to the inverse problem for the unperturbed data: from left to
right $z$= 1.035, 1.094, 1.191, 1.582.}}
  \label{fig2}
\end{figure}
\begin{figure}[h]
  \centering
%%%%%%%%%%%%%%%%%%%%%%%%%%%%%%%%%%%%%%%%%%%%%
\includegraphics[width=35mm,height=55mm]{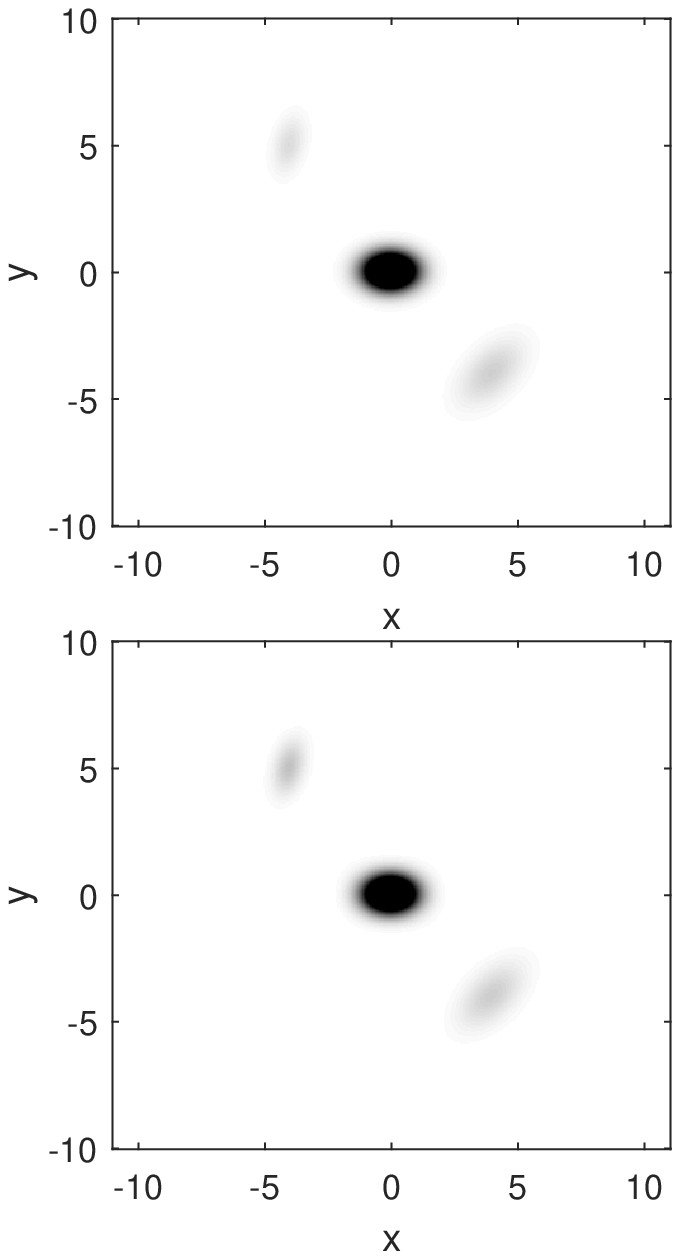}% bmp
\includegraphics[width=35mm,height=55mm]{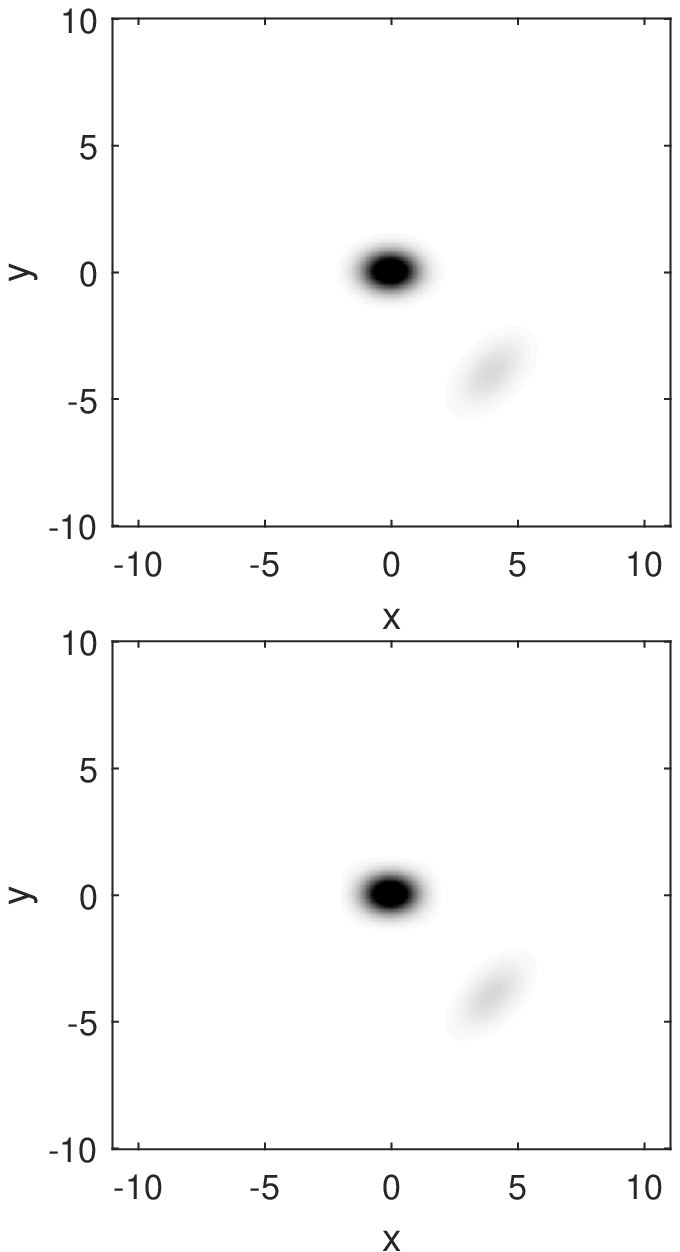}
\includegraphics[width=35mm,height=55mm]{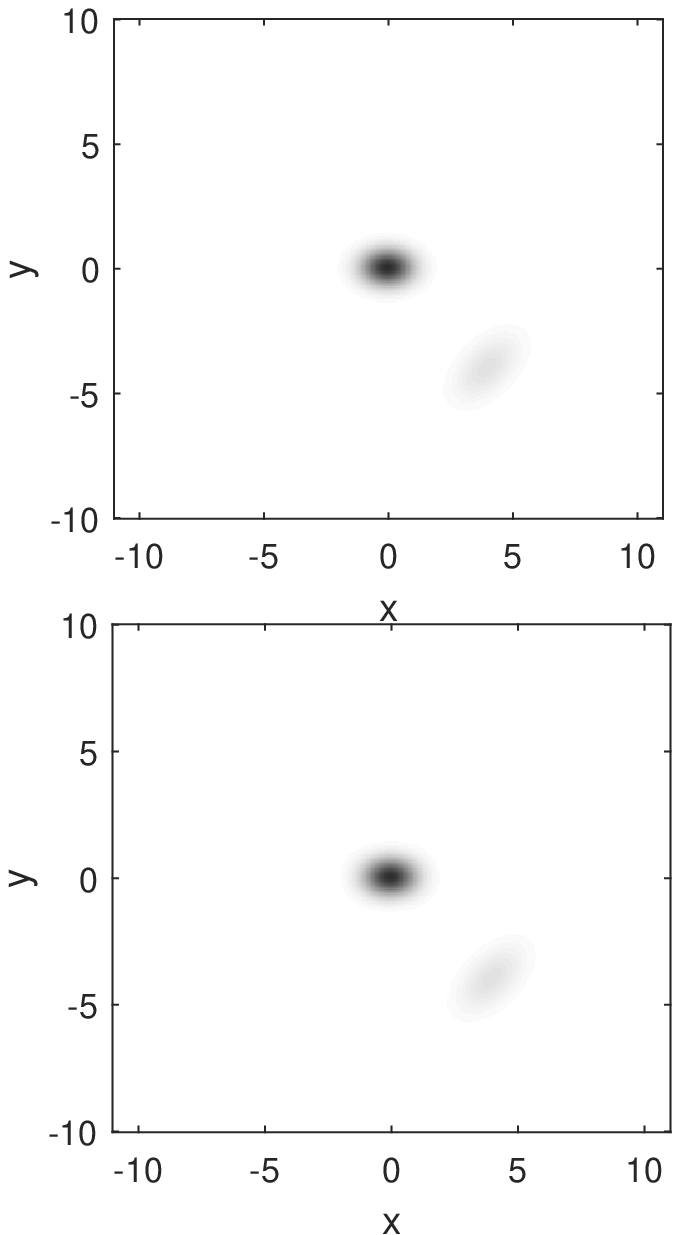}
\includegraphics[width=35mm,height=55mm]{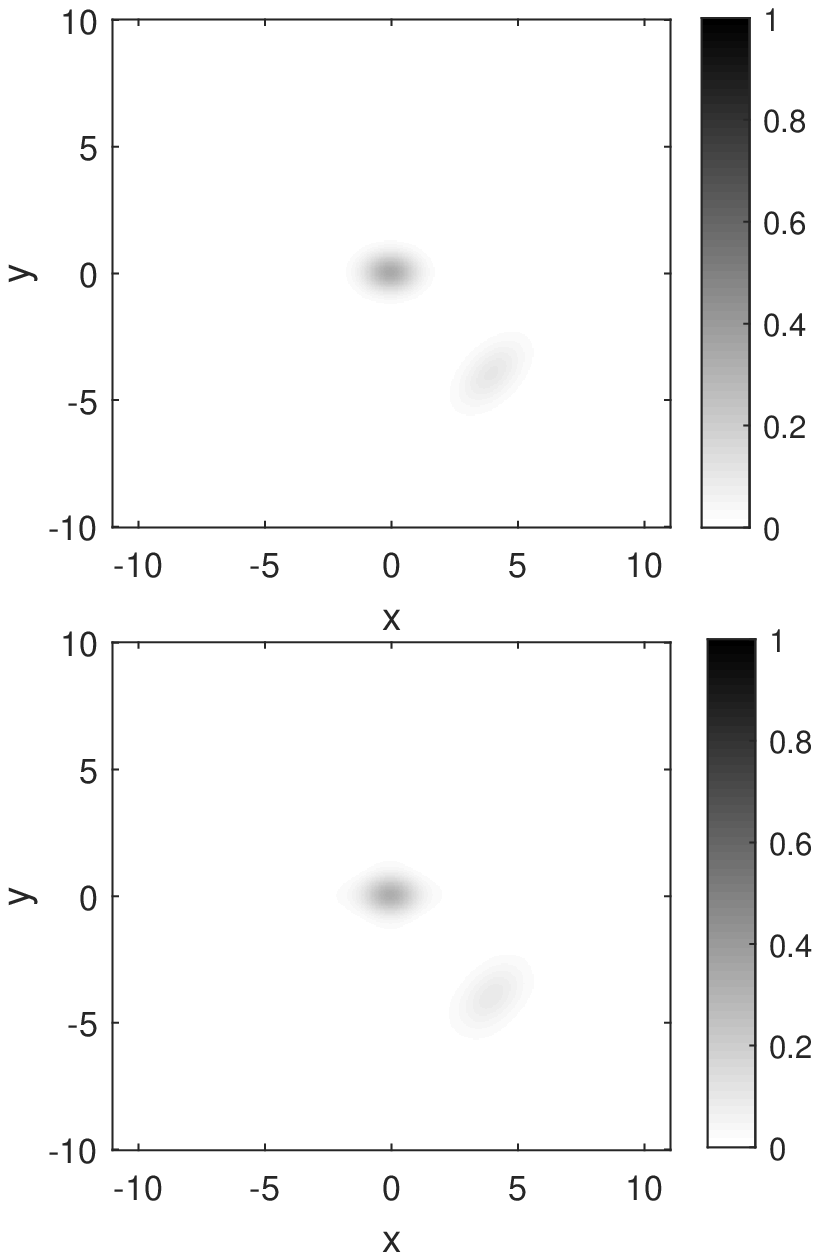}
%%%%%%%%%%%%%%%%%%%%%%%%%%%%%%%%%%%%%%%%%%%%%
  \caption{{\small A qualitative comparison of the exact (top) and an approximate
(bottom) solutions to the inverse problem for the unperturbed data: from left to
right $z$= 1.777, 1.875,
1.934, 1.973.}}
  \label{fig3}
\end{figure}
We also assume that $g(t)=\exp (-t)$ and $c_0=0.5$. Then we
compute the discrete analogue of the experimental data
\begin{multline*}
\frac{\partial V_{2}(x)}{\partial x_{3}}=\frac{1}{2\pi }\{ \int_{R}%
\frac{\partial }{\partial x_{3}}\left( \frac{1}{\left\vert
x-x^{\prime
}\right\vert }\right) \zeta (x^{\prime })dx^{\prime }-\frac{A_{2}}{2}%
\sum_{l=1}^{L}\frac{\partial }{\partial x_{3}}\left(
\frac{1}{\left\vert
x-x_{l}\right\vert }\right) +\\+\frac{A_{0}}{4\pi c_{0}^{2}}\sum_{l=1}^{L}%
\int_{\Pi }\frac{\partial }{\partial x_{3}}\left(
\frac{1}{\left\vert x-x^{\prime }\right\vert }\right)
\frac{dx^{\prime }}{\left\vert x^{\prime }-x_{l}\right\vert
}\}, \, \, \, x\in Y,
\end{multline*}
in (\ref{EQ7}), using a discrete Fourier transform
(algorithmically, fast Fourier transform) as described above.
In so doing, we obtain the data of the inverse problem, $\frac{%
\partial V_{2}(x)}{\partial x_{3}}$, with errors associated with finite-dimensional
approximation in calculation of the integrals. After that we
form the right-hand sides of the discrete equations (\ref{EQ6})
and solve them by the TSVD method. Finally, we calculate the
function $\xi (x,y,z)=\frac{\zeta(x,y,z)}{V_0(x,y,z)}$. We
confine ourselves to the calculation of this particular
function, not $ c(x) $, because examine the procedure for
solving linear equations (\ref{EQ7}).

Figures 2 and 3 show, for several values of $z\in [1,2]$, the
exact solution $\xi _{exact} (x,y,z)$ (top) and calculated
approximate solution $\xi _{appr} (x,y,z)$ (below) in pairs for
qualitative comparisons. The quantitative dependence of
obtained relative accuracy,
\[\Delta _{C} (z)=\frac{\left\| \xi _{appr} (x,y,z)-\xi _{exact}
(x,y,z)\right\| _{C(\mathbb{R}_{xy}^{2} )} }{\left\| \xi _{exact}
 (x,y,z)\right\| _{C(\mathbb{R}_{xy}^{2} )} }, \]
for the approximate solution $\xi _{appr} (x,y,z)$ in the layer
$z=\mathrm{const}$ is shown in Fig. 4 by solid line. The
calculations were performed in MATLAB on PC with a processor
Intel (R) Core (TM) i7-2600 CPU 3.40 GHz, RAM 8GB without
parallelization.
\begin{figure}[h]
  \centering
%%%%%%%%%%%%%%%%%%%%%%%%%%%%%%%%%%%%%%%%%%%%%
\includegraphics[width=70mm,height=50mm]{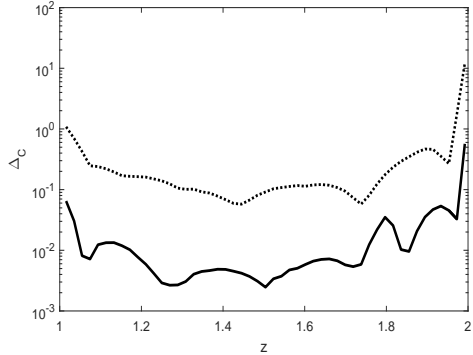}% bmp
%%%%%%%%%%%%%%%%%%%%%%%%%%%%%%%%%%%%%%%%%%%%%
  \caption{{\small The relative accuracy $\Delta _{C} (z)$ of approximate
solution (in the norm of $C(\mathbb{R}_{xy}^{2} )$) for different $z$. Solid
line: the value $\Delta _{C} (z)$ for for the unperturbed data, dotted line:
$\Delta _{C} (z)$ for perturbed data with a maximum perturbation value
about 1e-8.}}
  \label{fig5}
\end{figure}
\begin{figure}[h]
  \centering
%%%%%%%%%%%%%%%%%%%%%%%%%%%%%%%%%%%%%%%%%%%%%
\includegraphics[width=70mm,height=50mm]{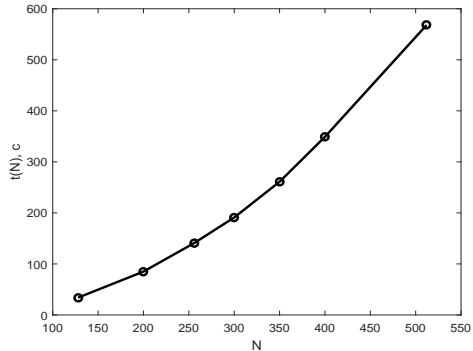}% bmp
%%%%%%%%%%%%%%%%%%%%%%%%%%%%%%%%%%%%%%%%%%%%%
  \caption{{\small The solving time vs $N$.}}
  \label{fig5}
\end{figure}

Since we are talking about the creation of a fast solution
algorithm for our inverse problem, we point out some of its
temporal characteristics. The dimensions of grids in the
variables $z,z'$, i.e, the values of $M,\,M'$, determine the
speed of solving the one-dimensional integral equation
(\ref{EQ5}) at fixed $\omega _{1} ,\omega _{2} $, that is the
SLAE (\ref{EQ6}) for a fixed $m$. The corresponding solution
time, $t_{0} (M,M')$, varies a little in passing from one of
equations (\ref{EQ5}) to another. This time is controlled by
the desired resolution of the algorithm in the variable $z$.
Actually, we have the estimate $t(N,M,M')\approx t_{0}
(M,M')\cdot N^{2} $ for full-time inverse problem solution on
chosen grids, and the number $N$ is controlled here by required
resolution in $x,y$. We present in Figure 5 the dependence
$t(N)=t(N,M,M')$ calculated on the same computer for fixed
$M=51,M'=51$ and for different $N$. In particular, the time to
solve the 3D inverse problem for $N=512$ is less than 10
minutes.

Note that the inverse problem under consideration is extremely
sensitive to errors in the input data. In solving it with the
double precision, small perturbations in the right side of the
SLAEs (\ref{EQ6}) by random errors with the level of the order
of 1e-8 lead to excessive smoothing of the approximate solution
when used the TSVD method and Tikhonov regularization.

The reason is very fast decay of singular values of the
matrices $A^{(m)}$ in the system (\ref{EQ6}). The graph
characterizing their behavior for a typical matrix $A^{(m)}$ is
shown in Figure 6. When using the TSVD method for the solution
of SLAEs (\ref{EQ6}), the small singular values are rejected
\cite{Hans,EHN}. In our case, for the inverse problem data
computed with the approximation error, but without the
additional perturbation by random error, singular values of the
order 1e-12 -- 1e-13 or less are discarded. The remainder of a
singular basis, 15 - 25 elements are used in the solution of
(\ref{EQ6}) and allow to reproduce the desired solution of the
inverse problem relatively accurate (see Figure 4, solid line).
The introduction of random errors of the order of 1e-10 -- 1e-8
in the data drastically reduces the dimension of used singular
basis (to 4 -- 5 elements), and this leads to a
''oversmoothed'' approximate solution, i.e. to its poor
accuracy (see Figure 4, dotted line). Approximately the same
effect occurs when using Tikhonov regularization. However, the
degree of oversmoothness for approximate solution is more than
for the TSVD. The corresponding theoretical error estimates
under different a priori assumptions can be found in
\cite{2,3}.
\begin{figure}[h]
  \centering
%%%%%%%%%%%%%%%%%%%%%%%%%%%%%%%%%%%%%%%%%%%%%
\includegraphics[width=70mm,height=60mm]{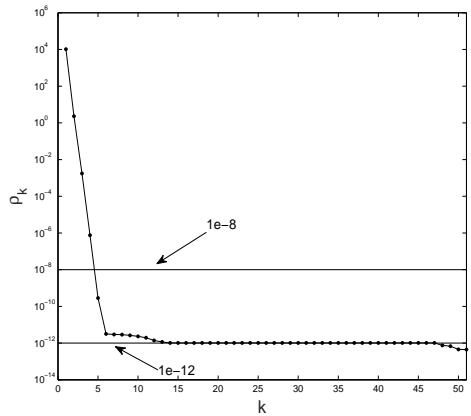}% bmp
%%%%%%%%%%%%%%%%%%%%%%%%%%%%%%%%%%%%%%%%%%%%%
  \caption{{\small Dependence of singular values $\rho _{k} $ of a matrix $A^{(m)} $
  on their number $k$.}}

  \label{fig6}
\end{figure}

Summing up the results of this work, we can draw the following
conclusions.

1. The inverse coefficient problem for the wave equation,
arising in modeling acoustic sensing, can be solved numerically
faster if instead of the time dependence $u(x,t)$ of the
scattered field in a region $Y$ we take as input data some
integrals of $u(x,t)$ in time. One possible integral of such a
kind is the function $V_{2} (x)$ or its partial derivatives.

2. For this type of data recorded in the plane layer, it is
possible to propose a numerical algorithm, which allows to
solve the inverse problem on a personal computer without the
use of supercomputer systems, in a relatively short period of
time, for sufficiently fine grids.

The work was supported by the Russian Foundation for Basic
Research (grants nos. 16-01-00039 and 15-01-00026а).

\end{document}